\documentclass[a4paper]{article}

\usepackage{amsmath}
\usepackage{amsfonts}
\usepackage{amssymb}
\usepackage{amsthm}
\usepackage{epsfig}
\usepackage{url}

\renewcommand{\emph}[1]{\textbf{{#1}}}
\def\endrem{}

\newcommand{\poly}{\mathcal{P}}
\newcommand{\NP}{\mathcal{NP}}
\newcommand{\coNP}{\text{co-}\mathcal{NP}}
\newcommand{\kverts}[2]{\mathcal{V}_{#1}\left({#2}\right)}

\newcommand{\Hsum}[1]{\mathcal{H}\left({#1}\right)}
\newcommand{\Hksum}[2]{\mathcal{H}_{#1}\left({#2}\right)}

\newcommand{\graph}[1]{\mathcal{G}_{{#1}}}
\newcommand{\verts}[1]{\mathcal{V}\left({#1}\right)}
\newcommand{\edges}[1]{\mathcal{E}\left({#1}\right)}
\newcommand{\dual}[1]{{#1}^{\Delta}}

\newtheorem{thm}{Theorem}
\newtheorem{lem}[thm]{Lemma}
\newtheorem{cor}[thm]{Corollary}
\newtheorem{defn}[]{Definition}

\newcommand{\calO}{\mathcal{O}}
\newcommand{\calS}{\mathcal{S}}
\newcommand{\prob}{\mathcal{D}}

\title{On the $k$-Systems of a Simple Polytope}

\author{ Michael Joswig\thanks{Partially supported by Deutsche
    Forschungsgemeinschaft, Sonderforschungsbereich~288
    ``Differentialgeometrie und Quantenphysik''.}  
  \and Volker Kaibel\thanks{Supported by Deutsche
    Forschungsgemeinschaft, Gerhard-Hess-Forschungsf\"orderungspreis
    G\"unter M.~Ziegler (Zi 475/2-3).}  
  \and
  Friederike K\"orner }

\begin{document}

\maketitle

\begin{abstract}
  A $k$-system of the graph~$\graph{P}$ of a simple polytope~$P$ is a
  set of induced subgraphs of~$\graph{P}$ that shares certain
  properties with the set of subgraphs induced by the $k$-faces
  of~$P$. 
  This new concept leads to polynomial-size certificates in
  terms of~$\graph{P}$ for both the set of vertex sets of facets as
  well as for abstract objective functions (AOF) in the sense of
  Kalai. 
  Moreover, it is proved that an acyclic orientation yields an
  AOF if and only if it induces a unique sink on every $2$-face.
\end{abstract}


\vspace{.3cm}

\noindent
\textbf{MSC~2000:}
  52B11 ($n$ dimensional polytopes), 52B05 (combinatorial properties)


\section{Introduction}
\label{sec:intro}

A celebrated theorem of Blind and Mani~\cite{BM87} states that the
combinatorial type of any simple polytope~$P$ is determined by the
isomorphism class of its abstract vertex-edge graph~$\graph{P}$.
Kalai~\cite{Kal88} gave a short and very elegant proof of this result.
The proof is constructive, but the algorithm that can be derived from
it has a worst-case running time which is exponential in the size
of~$\graph{P}$ (for computational experiments see Achatz and
Kleinschmidt~\cite{AK00}).  Thus, the complexity status of the problem
of reconstructing the combinatorial type of a simple polytope from its
graph remains unclear.

Kalai's proof is based on an ingenious characterization of the
shellings of the boundary of the dual polytope~$\dual{P}$.  Each
shelling order of the facets of $\partial\dual{P}$ corresponds to a
linear extension of an acyclic orientation of~$\graph{P}$ which
induces a unique sink in each non-empty face.  Such a linear ordering
of the vertices is called an \emph{abstract objective function}, while
the corresponding orientation is an \emph{AOF-orientation}.  Abstract
objective functions generalize linear objective functions in general
position.  The crucial step in Kalai's proof is the characterization
of AOF-orientations as those acyclic orientations of~$\graph{P}$ which
minimize a certain integer-valued function~$\Hsum{\calO}$.  Its
minimum value is the total number of non-empty faces of~$P$.

We consider a refinement of the function
$\Hsum{\calO}=\sum_{k=0}^d\Hksum{k}{\calO}$ as a sum of~$(d+1)$
functions $\Hksum{0}{\calO},\dots,\Hksum{d}{\calO}$.  This refinement
becomes useful in connection with the concept of a $k$-system that we
propose. A $k$-system of the graph~$\graph{P}$ of a simple
$d$-polytope~$P$ is a set of induced subgraphs of~$\graph{P}$
satisfying simple combinatorial conditions (that can be checked in
polynomial time) that, in particular, are fulfilled by the set of
subgraphs induced by the $k$-faces. Our main result on $k$-systems
(Theorem~\ref{thm:Hksum}) is that on the one hand, a $k$-system of the
graph of~$P$ with maximal cardinality is the set of subgraphs induced
by the $k$-faces of~$P$, and on the other hand, there is a strong dual
relation between the cardinality of $k$-systems and the
function~$\Hksum{k}{\calO}$. From this relationship polynomially sized
proofs (certificates) for the fact that a set of vertex sets indeed is
the set of vertex sets of the $k$-faces are readily obtained.  Note
that these certificates are purely combinatorial. In particular, no
coordinates are involved.

Furthermore, we prove that every acyclic orientation which induces a
unique sink in every $2$-face of~$P$ is an AOF-orientation
(Theorem~\ref{thm:2faces}). This reveals a strong relationship between
the $2$-faces and the abstract objective functions of a simple
polytope; they can be exploited as certificates for each other.  The
special role which is played by the $2$-skeleton  reflects the
well-known fact that it is straightforward to reconstruct a simple
polytope from its $2$-skeleton.

We refer to Ziegler's book~\cite{Zie95} for a detailed treatment of
all notions and concepts we rely on.


\section{Results}
\label{sec:results}

Let~$P$ be a simple $d$-polytope. We denote the graph of~$P$
by~$\graph{P}=(\verts{P},\edges{P})$, where~$\verts{P}$ is the set of
vertices of~$P$ and~$\edges{P}$ is the set of its edges.

If~$W\subseteq\verts{P}$ is a subset of vertices, then~$\graph{P}(W)$
is the subgraph of~$\graph{P}$ induced by~$W$.  For each $0\leq k\leq
d-1$ let~$\kverts{k}{P}$ be the set of vertex sets of $k$-faces
of~$P$.  As usual, $f_k(P):=|\kverts{k}{P}|$ is the number of
$k$-faces of~$P$.  We will often identify a face~$F$ of~$P$ with the
subgraph of~$\graph{P}$ (denoted by ~$\graph{P}(F)$) that is induced
by the vertices of~$F$.

\begin{defn}
  Let~$P$ be a simple $d$-polytope and $2\leq k\leq d-1$.
  \begin{enumerate}
  \item[(i)] A \emph{\boldmath{k}-frame} of~$P$ is a (not necessarily
    induced) subgraph of~$\graph{P}$ isomorphic to the star~$K_{1,k}$,
    where the vertex of degree~$k\geq 2$ is called the \emph{root} of
    the $k$-frame.
  \item[(ii)] A set~$\calS$ of subsets of~$\verts{P}$ is called a
    \emph{\boldmath{k}-system} of~$\graph{P}$ if for every set $S\in\calS$ the
    subgraph~$\graph{P}(S)$ of~$\graph{P}$ is $k$-regular and the node
    set of every $k$-frame of~$P$ is contained in a unique set
    from~$\calS$.
  \end{enumerate}
\endrem
\end{defn}

\begin{figure}[ht]
  \begin{center}
    \epsfig{figure=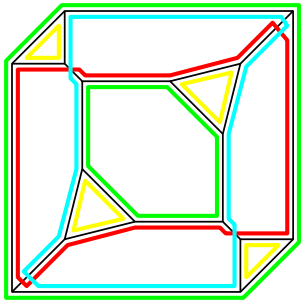,height=4cm}
    \caption{A $2$-system (indicated by the subgraphs induced by its
      sets) that is not the set of vertex sets of $2$-faces. The
      polytope arises from cutting off two opposite vertices of a
      $2$-face of the $3$-cube.}
    \label{fig:2system}
  \end{center}
\end{figure}

Obviously, $\kverts{k}{P}$ is a $k$-system.  In general,
$\kverts{k}{P}$ is not the only $k$-system of~$\graph{P}$.
Figure~\ref{fig:2system} shows a $2$-system of the graph of a simple
$3$-polytope~$P$ that is different from~$\kverts{2}{P}$. We will
characterize~$\kverts{k}{P}$ among the $k$-systems of~$\graph{P}$ by
means of certain acyclic orientations.

Let~$\calO$ be an acyclic orientation of~$\graph{P}$.  It is an
elementary (but crucial) fact that for every~$W\subseteq\verts{P}$ the
orientation of~$\graph{P}(W)$ induced by~$\calO$ has at least one
sink; furthermore, from each $w\in W$ there is a directed path in
$\graph{P}(W)$ to one of these sinks.  For every~$0\leq i\leq d$, let
$h_i(\calO)$ be the number of vertices of~$\graph{P}$ with
precisely~$i$ of its incident edges directed inwards. We define (for
all $0\leq k\leq d$)
$$
\Hksum{k}{\calO}:=\sum_{i=0}^d h_i(\calO)\binom{i}{k}
\qquad\text{and}\qquad \Hsum{\calO}:=\sum_{i=0}^d
h_i(\calO)2^i=\sum_{k=0}^d\Hksum{k}{\calO}\enspace.
$$

The sum~$\Hksum{k}{\mathcal{O}}$ is the number of $k$-frames for which
all edges are directed towards the root. Thus~$\Hksum{k}{\calO}$ is
the total number of sinks induced in the subgraphs~$\graph{P}(S)$
of~$\graph{P}$ ($S\in\calS$).

One of the
beautiful steps on Kalai's ``Simple Way to Tell a Simple Polytope from
its Graph''~\cite{Kal88} (see also~\cite{Zie95}, Chap.~3.4) is the
observation that the AOF-orientations of~$\graph{P}$ are precisely
those orientations that minimize~$\Hsum{\calO}$.
Theorem~\ref{thm:2faces} implies that AOF-orientations can also be
characterized as those acyclic orientations of~$\graph{P}$ that
minimize~$\Hksum{2}{\calO}$.  If~$\calO$ is an AOF-orientation
of~$\graph{P}$, then $(h_0(\calO),\dots,h_d(\calO))$ is the
$h$-vector of~$P$ (see, e.g.,~\cite{Zie95}, Chap.~8.3); in
particular, the numbers~$h_k(\calO)$ do not depend on the specific
choice of the AOF-orientation~$\calO$.

There is an important relationship between the $k$-systems
and the acyclic orientations of~$\graph{P}$.

\begin{thm}
  \label{thm:Hksum}
  Let~$P$ be a simple $d$-polytope, and let $2\leq k\leq d-1$. For
  every $k$-system~$\calS$ of~$\graph{P}$ and every acyclic
  orientation~$\calO$ of~$\graph{P}$ the inequalities
  $$
  \left|\calS\right| \stackrel{(1)}{\leq} f_k(P)
  \stackrel{(2)}{\leq} \Hksum{k}{\calO}
  $$
  hold, where~\mbox{\rm (1)} holds with equality if and only if
  $\calS=\kverts{k}{P}$, and~\mbox{\rm (2)} holds with equality if and only
  if~$\calO$ induces precisely one sink on every $k$-face of~$P$.
  \end{thm}

\textit{Proof.}
  Let~$\calS$ be a $k$-system of~$\graph{P}$, and let~$\calO$ be an
  acyclic orientation of~$\graph{P}$. Since~$\calO$ is acyclic, $\calO$
  induces at least one sink in every $S\in\calS$. In particular,
  $\Hksum{k}{O}\geq|\calS|$ holds.  Hence, inequality~(2) (together
  with the characterization of equality) follows with
  $\calS:=\kverts{k}{P}$, and inequality~(1) is obtained by
  choosing~$\calO$ as any AOF-orientation of~$\graph{P}$.  It remains
  to show that $|\calS|=f_k(P)$ implies $\calS=\kverts{k}{P}$.
  
  Let~$\calS$ be a $k$-system of~$\graph{P}$ with~$|\calS|=f_k(P)$. In order
  to show $\calS=\kverts{k}{P}$ it  suffices to prove
  $\kverts{k}{P}\subseteq\calS$. The main ideas of the following are
  imported from Kalai's paper~\cite{Kal88}.  Let~$W\in\kverts{k}{P}$
  be the vertex set of any~$k$-face~$F$ of~$P$. There is a linear
  function (in general position) on~$P$ which assigns larger values to
  the vertices on~$F$ than to all other vertices of~$P$. This linear
  function induces an AOF-orientation~$\calO$ of~$\graph{P}$ with the
  property that no edge is directed into~$W$ ($W$ is initial).
  
\begin{figure}[ht]
    \begin{center}
      \epsfig{figure=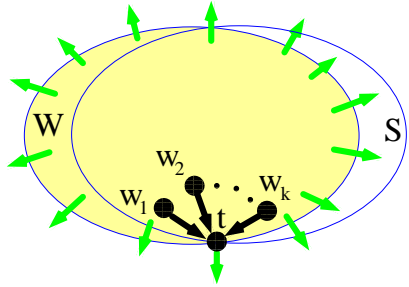, height=4cm}
      \caption{Illustration of the proof of Theorem~\ref{thm:Hksum}.}
      \label{fig:Hksum}
    \end{center}
  \end{figure}

  See Fig.~\ref{fig:Hksum} for a sketch of the situation.
  Denote by~$t\in W$ the unique sink induced by~$\calO$
  in~$\graph{P}(W)$, and let $w_1,\dots,w_k\in W$ be the neighbors
  of~$t$ in~$F$. Let~$S$ be the (unique) set in~$\calS$ containing the
  $k$-frame with node set~$\{t,w_1,\dots,w_k\}$. Due to
  $|\calS|=f_k(P)=\Hksum{k}{\calO}$ the orientation~$\calO$ induces a
  unique sink in~$\graph{P}(S)$, which must be~$t$.
  Since~$W$ is initial, this implies~$S\subseteq W$ (because there
  must be a directed path from every vertex in~$S$ to~$t$).  Hence,
  $\graph{P}(S)$ is a $k$-regular subgraph of the $k$-regular and
  connected graph~$\graph{P}(W)$.  Thus, $W=S\in\calS$.
\qed

As a consequence of Theorem~\ref{thm:Hksum}, we obtain a
characterization of~$\kverts{k}{P}$ that is quite similar to Kalai's
characterization of AOF-orientations via minimizers
of~$\Hsum{\calO}$. 

\begin{cor}
  \label{cor:maxkSystem}
  Let~$P$ be a simple $d$-polytope. A  $k$-system of~$\graph{P}$
  is~$\kverts{k}{P}$ if and only if it has maximal cardinality among
  all $k$-systems of~$\graph{P}$.
\end{cor}

Similar to Kalai's result, this corollary implies that the
set~$\kverts{d-1}{P}$ of vertex sets of facets of~$P$ can be computed
from~$\graph{P}$.  However, it does not shed any light on the
question, how fast this can be done. From the complexity point
of view, the next characterization (which follows from
Theorem~\ref{thm:Hksum}, since every simple $d$-polytope has an
AOF-orientation) is much more valuable.

\begin{cor}
  \label{cor:Hksum}
  Let~$P$ be a simple $d$-polytope, and let~$\calS$ be a $k$-system
  of~$\graph{P}$ (with $2\leq k\leq d-1$). Then either there is an
  acyclic orientation~$\calO$ of~$\graph{P}$ with
  $\Hksum{k}{\calO}=|\calS|$ or there is a $k$-system~$\calS'$
  of~$\graph{P}$ with~$|\calS'|>|\calS|$. In the first case,
  $\calS=\kverts{k}{P}$, in the second, $\calS\not=\kverts{k}{P}$.
\end{cor}

Corollary~\ref{cor:Hksum} yields a good characterization
of~$\kverts{k}{P}$ among all sets of subsets of vertices of a simple
polytope~$P$ (given by its graph) in the sense of
Edmonds~\cite{Edm65a,Edm65b}: for every subset~$\calS$ of vertex sets
of~$P$ one can efficiently prove the answer to the question
``Is~$\calS=\kverts{k}{P}$?'' (although it is currently unknown if one
can also find the answer efficiently). If the answer is ``yes,'' then
we may prove this in polynomially many steps (in the size
of~$\graph{P}$) by first checking that~$\calS$ is a $k$-system, and
then exhibiting an acyclic orientation~$\calO$ of~$\graph{P}$ with
$|\calS|=\Hksum{k}{\calO}$. If the answer is ``no,'' then we may prove
this by showing that~$\calS$ is not a $k$-system of~$\graph{P}$, or,
if it is a $k$-system, by exhibiting a larger $k$-system~$\calS'$
of~$\graph{P}$.

Since the number of facets of a simple $d$-polytope~$P$ is bounded by
a polynomial in the size of~$\graph{P}$, Corollary~\ref{cor:Hksum}
also implies that the question whether a given single subset of
vertices is the vertex set of some facet of~$P$ has a good
characterization.

It had been hoped for a long time that such good characterizations
(for $k=d-1$) would be obtained by an eventual proof of a conjecture
due to Perles.  Already in~1970 he conjectured that every subset of
vertices of a simple $d$-polytope~$P$ which induces a $(d-1)$-regular,
connected, non-separating subgraph of~$\graph{P}$ is the vertex set of
a facet of~$P$. This would even imply much more than good
characterizations: it would immediately yield polynomial time
algorithms to decide whether a set of vertices is the vertex
set of a facet, and whether a subset of sets of vertices
is~$\kverts{d-1}{P}$. However, recently Haase and Ziegler~\cite{HZ00}
disproved Perles' conjecture.

For $k=2$, Theorem~\ref{thm:Hksum} (together with
Theorem~\ref{thm:2faces}) also provides us with a good
characterization of the AOF-orientations among all acyclic
orientations of the graph~$\graph{P}$ of a simple polytope~$P$ (see
Corollary~\ref{cor:H2sum}).  Previously, the only method that was
known to prove that an acyclic orientation of~$\graph{P}$ is an
AOF-orientation was to show that it minimizes~$\Hsum{\calO}$ by
exploring all acyclic orientations of~$\graph{P}$ (where it was
perhaps the most striking result of~\cite{Kal88} that such a method
does exist at all). Notice that, if in addition to~$\graph{P}$
also~$\kverts{d-1}{P}$ is specified as input, then it can be decided
in polynomial time whether an acyclic orientation of~$\graph{P}$ is an
AOF-orientation. This follows easily from the equivalence between
AOF's on~$P$ and shellings of~$\dual{P}$.  However, in our context the
polytope~$P$ is specified just by its graph, and the ultimate question
is whether~$\kverts{d-1}{P}$ can be computed efficiently at all.

\begin{cor}
  \label{cor:H2sum}
  Let~$P$ be a simple $d$-polytope, and let $\calO$ be an acyclic
  orientation of~$\graph{P}$. Then either there is a
  $2$-system~$\calS$ of~$\graph{P}$ with $|\calS|=\Hksum{2}{\calO}$ or
  there is an acyclic orientation~$\calO'$ of~$\graph{P}$ with
  $\Hksum{2}{\calO'}<\Hksum{2}{\calO}$.  In the first case, $\calO$ is
  an AOF-orientation, in the second, it is not.
\end{cor}

While the ``either or''-statement of the Corollary follows immediately
from Theorem~\ref{thm:Hksum}, the fact that in the first case $\calO$
is an AOF-orientation is implied by the following result (which, in
particular, implies that Ex.~8.12~(iv) in~\cite{Zie95} cannot be
solved). The theorem had already been proved for hypercubes
by~\cite{HSLdW88}. For $3$-dimensional simple polytopes the result of
Theorem~\ref{thm:2faces} has independently been obtained by
Develin~\cite{Dev00}.

\begin{thm}
  \label{thm:2faces}
  Let~$P$ be a simple polytope, and let~$\calO$ be an acyclic
  orientation of~$\graph{P}$. If~$\calO$ induces precisely one sink on 
  every $2$-face of~$P$, then~$\calO$ is an AOF-orientation.
\end{thm}

Since every face of a simple polytope is simple,
Theorem~\ref{thm:2faces} follows immediately from the following
result.

\begin{lem}
  \label{lem:sinksInFacets}
  Let~$P$ be a simple~$d$-polytope and~$2\leq k\leq d-1$. If~$\calO$
  is an acyclic orientation of~$\graph{P}$ that has more than one
  global sink, then there is a $k$-face of~$P$ on which~$\calO$
  induces more than one sink.
\end{lem}

\begin{figure}[ht]
    \begin{center}
      \epsfig{figure=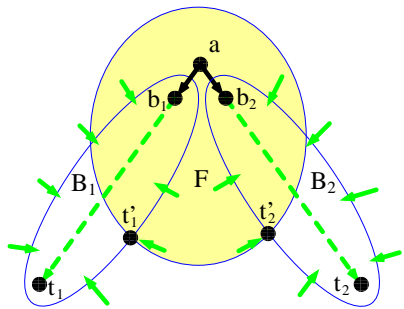, height=4cm}
      \caption{Illustration of the proof of Lemma~\ref{lem:sinksInFacets}.}
      \label{fig:sinksInFacets}
    \end{center}
  \end{figure}

\textit{Proof.}
  Let $\calO$ have more than one sink in~$\graph{P}$.  We denote by
  $A\subseteq\verts{P}$ the set of all vertices from which two
  different sinks can be reached on directed paths.  Since~$\graph{P}$
  is connected, $A\not=\varnothing$. Thus we can choose a vertex $a\in
  A$ which is a sink in~$\graph{P}(A)$, together with two directed
  paths $(a,b_1,\dots,t_1)$ and $(a,b_2,\dots,t_2)$ connecting~$a$
  with two distinct (global) sinks~$t_1$ and~$t_2$ (see
  Fig.~\ref{fig:sinksInFacets} for an illustration of the proof).
  
  Since~$P$ is simple and $k\geq 2$, there is a $k$-face~$F$
  containing $a$, $b_1$, and $b_2$.  For $i\in\{1,2\}$ denote by~$B_i$
  the set of vertices of~$\graph{P}$ that lie on some directed path
  from~$b_i$ to~$t_i$.  The choice of~$a$ as a sink in~$\graph{P}(A)$
  implies $B_1\cap B_2=\varnothing$.  Since both $B_1\cap F$ and
  $B_2\cap F$ are non-empty, the acyclic orientation~$\calO$ thus
  induces two distinct sinks~$t'_1$ and~$t'_2$
  in~$\graph{P}(B_1\cap F)$ and~$\graph{P}(B_2\cap F)$, respectively.
  Again, since~$a$ is a sink in~$\graph{P}(A)$, both $B_1$ and~$B_2$
  are terminal (no edges are directed outwards).  Hence $t'_1$
  and~$t'_2$ are two distinct sinks in~$\graph{P}(F)$ as well.
\qed

It is not too hard to find examples showing that there is no analogue
of Theorem~\ref{thm:2faces} for~$k$-faces with $k>2$.
Theorem~\ref{thm:2faces} thus shows that the $2$-faces of a simple
polytope~$P$ play a distinguished role with respect to the
AOF-orientations of~$\graph{P}$. It is worth noticing that
computing~$\kverts{2}{P}$ from~$\graph{P}$ is polynomial time
equivalent to computing~$\kverts{d-1}{P}$ from~$\graph{P}$. To see
this, it suffices to observe that the obvious bijections between the
neighbors of~$u$ and the neighbors of~$v$ (for every edge
$\{u,v\}\in\edges{P}$) defined by the $2$-faces and the facets of~$P$,
respectively, coincide. Thus, instead of considering the problem of
computing~$\kverts{d-1}{P}$ from~$\graph{P}$ one may rather consider
the problem of computing~$\kverts{2}{P}$ from~$\graph{P}$.  The
$2$-faces are polygons and thus have a simpler structure than the
facets, in general. Moreover, they also bear strong connections to the
AOF-orientations as stated in Theorem~\ref{thm:2faces}.


\section{Discussion}
\label{sec:disc}

A good characterization, as provided by Corollary~\ref{cor:Hksum},
often indicates that the corresponding (decision) problem~$\prob$
(given the graph~$\graph{P}$ of a simple polytope~$P$ and a
$k$-system~$\calS$ of~$\graph{P}$; is~$\calS=\kverts{k}{P}$?) can be
solved in polynomial time.  In fact, there are many examples of
combinatorial optimization problems, for which such a good
characterization has guided the algorithm design (primal-dual
algorithms). In the theory of computational complexity, this
corresponds to the fact that for most problems which are known to be
contained in the complexity class $\NP\cap\coNP$ it is even known that
they belong to the class~$\poly$ of problems solvable in polynomial
time (the most prominent exception is the problem of deciding whether
an integer number is a prime).

Unfortunately, Corollary~\ref{cor:Hksum} does not imply that
problem~$\prob$ is contained in $\NP\cap\coNP$, since it is unknown if
one can prove resp.\ disprove efficiently that a given graph is
(isomorphic to) the graph of some simple $d$-polytope. This question
is closely related to the Steinitz problem, the problem to decide
whether a given lattice is (isomorphic to) the face-lattice of some
polytope (with real-algebraic coordinates). The Steinitz problem is
known to be $\NP$-hard even in dimension four~(Theorem~9.1.2
in~\cite{RG96}). Furthermore, again already in dimension four, there
is no polynomial certificate for the Steinitz problem by specifying
coordinates~(Theorem~9.3.3 in~\cite{RG96}). This can be interpreted as
indications for the non-existence of a good characterization for the
``integrity'' of the input data of problem~$\prob$.  Thus, the good
characterization of~Corollary~\ref{cor:Hksum} seems to have no direct
complexity theoretical implications. Nevertheless, it might be
encouraging or even be exploited for the design of a polynomial
time algorithm for problem~$\prob$.

Theorem~\ref{thm:Hksum} shows that for solving problem~$\prob$ in
polynomial time it would suffice to design a polynomial time method
for computing~$f_k(P)$ from~$\graph{P}$. One way to achieve this could
be a polynomial time method for finding any AOF-orientation
of~$\graph{P}$.  However, it is not even known whether there is a
polynomial time algorithm for finding an AOF of a simple
$d$-polytope~$P$ given by its entire face-lattice (not even
for~$d=4$).  Equivalently, there is no polynomial time algorithm known
that finds a shelling of an abstract simplicial complex of
which one knows that it is isomorphic to the boundary complex
of a simplicial polytope.  Thus, an interesting question is the one
for alternative ways to calculate~$f_k(P)$ from~$\graph{P}$. For
instance, it might be easier to find a polynomial algorithm that finds
an acyclic orientation of~$\graph{P}$ which has only one sink per
$k$-face, from which one would obtain~$f_k(P)$ as well.

These considerations concern the problem of deciding whether a
set of candidates actually is the set~$\kverts{d-1}{P}$ of vertex sets
of facets (or, more generally, of $k$-faces) of a simple polytope~$P$
specified by its graph~$\graph{P}$. The genuine question, however, is
whether there is a polynomial time algorithm for
finding~$\kverts{d-1}{P}$ from~$\graph{P}$.
Corollary~\ref{cor:maxkSystem} shows that one can phrase this
problem as a maximization problem.  Hence, it might well be that
concepts and tools from Combinatorial Optimization (such as the
primal-dual method mentioned above) can help to eventually find a
``fast way to tell a simple polytope from its graph.''


\providecommand{\bysame}{\leavevmode\hbox to3em{\hrulefill}\thinspace}


\vfill
\noindent
Michael Joswig, Volker Kaibel, Friederike K\"orner\\
Technische Universit\"at Berlin\\
Fakult\"at~II, Institut f\"ur Mathematik\\
MA~6--2\\
Stra\ss e des 17.~Juni~136\\
10623~Berlin\\
Germany\\
\url{{joswig,kaibel,koerner}@math.tu-berlin.de}

\end{document}